\documentclass[aos,MSNbibl,nameyear,seceqn,dvips]{arximspdf}
\usepackage{dcolumn}
\usepackage{newsym2e}

\doi{10.1214/12-AOS987} 
\volume{40}
\issue{2}
\pubyear{2012}
\firstpage{891}
\lastpage{907}

\makeatletter
\newtheorem{theorem}{Theorem}[section]
\newtheorem{lemma}[theorem]{Lemma}
\newproclaim{definition}[theorem]{Definition} 
\newproclaim{example}[theorem]{Example}
\newcommand{\x}{\mathbf{x}}
\newcommand{\y}{\mathbf{y}}
\newcommand{\one}{\mathbf{1}}
\newcommand{\two}{\mathbf{2}}
\newcommand{\zero}{\mathbf{0}}
\newcolumntype{d}[1]{D{.}{.}{#1}}
\newcolumntype{e}[1]{D{.}{}{#1}}
\renewcommand{\pmod}[1]{\ (\operatorname{mod}#1)}
\makeatother

\begin{document}
\begin{frontmatter}

\title{Uniform fractional factorial designs}
\runtitle{Uniform fractional factorial designs}

\begin{aug}
\author[A]{\fnms{Yu} \snm{Tang}\thanksref{t2}\ead[label=e1]{ytang@suda.edu.cn}},
\author[B]{\fnms{Hongquan} \snm{Xu}\corref{}\thanksref{t3}\ead[label=e2]{hqxu@stat.ucla.edu}}
\and
\author[C]{\fnms{Dennis K. J.} \snm{Lin}\ead[label=e3]{dkl5@psu.edu}}
\runauthor{Y. Tang, H. Xu and D. K. J. Lin}
\thankstext{t2}{Supported by NNSF of China with Grant 10801104 and
China Scholarship Council.}
\thankstext{t3}{Supported by NSF Grants DMS-08-06137 and 11-06854.}
\affiliation{Soochow University,
University of California, Los Angeles
and\\ Pennsylvania State University}
\address[A]{Y. Tang\\School of Mathematical Sciences \\
Soochow University \\
Suzhou, Jiangsu 215006\\ China \\
\printead{e1}} 
\address[B]{H. Xu\\Department of Statistics \\
University of California \\
Los Angeles, California 90095-1554\\ USA \\
\printead{e2}}
\address[C]{D. K. J. Lin\\Department of Statistics \\
Pennsylvania State University \\
University Park,
Pennsylvania 16802\\ USA \\
\printead{e3}}
\end{aug}

\received{\smonth{6} \syear{2011}}
\revised{\smonth{2} \syear{2012}}

%
\begin{abstract}
The minimum aberration criterion has been frequently used in
the selection of fractional factorial designs with nominal factors.
For designs with quantitative factors, however, level permutation
of factors could alter their geometrical structures and statistical properties.
In this paper uniformity is used to further distinguish fractional
factorial designs, besides the minimum aberration criterion.
We show that minimum aberration designs have low discrepancies on average.
An efficient method for constructing uniform minimum aberration designs
is proposed and
optimal designs with 27 and 81 runs are obtained for practical use.
These designs have good uniformity {and are effective for studying
quantitative factors}.
\end{abstract}

%

\end{frontmatter}
%

\section{Introduction}
The minimum aberration criterion [\citet{FH1980}]
has been frequently used in the selection of regular fractional
factorial (FF) designs with nominal factors, as it provides nice design
properties.
This is especially important
when the experimenter has little knowledge about
the potential significance of factorial effects. {The readers are
referred to \citet{MW2006} and \citet{WH2009} for
existing theory and results on minimum
aberration designs.
\citet{DT1999}, \citet{TD1999} and \citet{XW2001} further
proposed generalized minimum aberration criteria for comparing
nonregular fractional factorial designs}.

\citet{CW2001}
and \citet{FM2001} found
that designs may have different geometrical structures and
statistical properties, even though they share the identical
word-length pattern.
In view of this, \citet{CY2004} {pointed out} 
that the distinction
in the analysis objective and strategy for
experiments with nominal or quantitative factors requires different
selection criteria and classification methods. For designs
with quantitative factors,
they {proposed to describe design properties using the geometrical structures.}
Two designs are said to be \textit{geometrically isomorphic}
if one can be obtained from the other by a permutation of factors
and/or reversing the level order of one or more factors.
For example, consider the two designs in Table \ref{Tab_TwoDes}.
Design $A$ is a regular $3^{3-1}$ FF design with
$F_3 = F_1 + F_2 \pmod{3}$,
and Design $B$ is formed by $F'_3 = F_1 + F_2 +2\pmod{3}$.
It is obvious that
these two designs are \textit{combinatorially isomorphic} to each other,
because one can be obtained from the other by permuting the levels in
the third column [i.e., map (0, 1, 2) to (2, 0, 1)].
However, they have different geometrical {structures}
and thus are geometrically nonisomorphic.
Design $B$ contains the center run with all ones,
while Design $A$ does not. If we reverse the level order
[i.e., map (0, 1, 2) to (2, 1, 0)] for all three columns,
Design $B$ is invariant while Design $A$ is not.
These two designs have different statistical properties due to
their different geometrical {structures.}

\begin{table}
\caption{Two combinatorially isomorphic designs with different
geometrical structures}\label{Tab_TwoDes}
\begin{tabular*}{\textwidth}{@{\extracolsep{\fill}}lccccc@{}}
\hline
\multicolumn{3}{@{}c}{\textbf{Design}  $\bolds{A}$} &\multicolumn{3}{c@{}}{\textbf{Design}  $\bolds{B}$}
\\[-5pt]
\multicolumn{3}{@{}c}{\hrulefill} &\multicolumn{3}{c@{}}{\hrulefill}\\
{$\bolds{F_1}$} & {$\bolds{F_2}$} & {$\bolds{F_3}$} & {$\bolds{F_1}$} & {$\bolds{F_2}$} & {$\bolds{F'_3}$} \\
\hline
0 & 0 & 0 & 0 & 0 & 2 \\
0 & 1 & 1 & 0 & 1 & 0 \\
0 & 2 & 2 & 0 & 2 & 1 \\
1 & 0 & 1 & 1 & 0 & 0 \\
1 & 1 & 2 & 1 & 1 & 1 \\
1 & 2 & 0 & 1 & 2 & 2 \\
2 & 0 & 2 & 2 & 0 & 1 \\
2 & 1 & 0 & 2 & 1 & 2 \\
2 & 2 & 1 & 2 & 2 & 0 \\
\hline
\end{tabular*}
\vspace*{3pt}
\end{table}

To further classify geometrically nonisomorphic designs,
\citet{CY2004} generalized the concept of
minimum aberration and used indicator function
to define the $\beta$-word-length pattern
based on a polynomial model.
Despite its theoretical beauty,
the complexity of indicator function
prohibits its use for design construction.
On the other hand, \citet{FM2001}
suggested using uniformity to compare the performance
of geometrically nonisomorphic designs.
Various discrepancies have been used as
measures of uniformity; {see \citet{FLWZ2000} and \citet{FLS2006}}. These discrepancies all have their
geometrical meanings and can be interpreted as the difference between
the empirical distribution and the uniform distribution. Among them,
the centered $L_2$-discrepancy (CD), proposed by
\citet{Hickernell1998}, is the most frequently used.

Both the $\beta$-word-length pattern and the centered $L_2$-discrepancy
reflect the geometrical structure of the design. Here we use the discrepancy
to choose FF designs mainly for two reasons. First, the centered
$L_2$-discrepancy
has a simple analytic formula; it is much faster to
calculate the discrepancy than the $\beta$-word-length pattern. The difference
between the computational times is substantial.
The second and more important reason is that the $\beta$-word-length pattern
is model-dependent while the centered $L_2$-discrepancy is model free.
\citet{CY2004} defined the $\beta$-word-length pattern based on a
polynomial model
but they further pointed out that the $\beta$-word-length pattern
needs to be {modified} in other situations.
Optimal designs constructed based on the $\beta$-word-length pattern
would rely on the specific model used.
In contrast, designs with low discrepancy tend to have good space
filling properties
and are model robust in the sense that they can guard against inaccurate
estimates caused by model misspecification [\citet{HL2002}].

Here we propose to construct uniform FF designs from
existing minimum aberration designs via level permutations.
Obviously, for two-level designs, there is no difference when levels
are permuted,
but for high-level designs, there are many unknowns to be studied.
For convenience, we will focus on three-level designs in this paper,
but the basic ideas can be extended to higher-level designs.

The paper is organized as follows.
In Section \ref{Sec_Relationship}, {we obtain a key theorem and show that
minimum aberration designs tend to have low discrepancy on average.
Then we introduce the concept of uniform minimum aberration design.
In Section~\ref{Sec_ThreeLevel}, we present an efficient way for
constructing three-level
regular uniform FF designs and construct uniform minimum aberration
designs with 27 runs and 81 runs for practical use.
These newly-constructed designs often outperform
existing uniform designs, especially when the number of factors is large.
In Section~\ref{Sec_Justify}, we examine the relationship between the
discrepancy and the
$\beta$-word-length pattern. Uniform minimum aberration designs appear
to perform well with respect to the $\beta$-word-length pattern.}
The last section gives a brief conclusion. For clarity, we defer all
proofs to the \hyperref[app]{Appendix}.

\section{Uniform minimum aberration designs}\label{Sec_Relationship}

A design with $N$ runs, $n$ factors and $s$ levels, denoted by $(N, s^n)$,
is an $N\times n$ matrix. Throughout the paper,
the $s$ levels are denoted as $0, 1, \ldots, s-1$.
For an $(N,s^n)$-design $D$, consider an ANOVA model
\[
Y = X_0 \alpha_0 + X_1 \alpha_1 + \cdots+ X_n \alpha_n + \varepsilon,
\]
where $Y$ is the vector of $N$ observations, $\alpha_0$ is the
intercept and $X_0$ is an $N\times1$ vector of $1$'s, $\alpha_j$ is
the vector of all $j$-factor interactions and $X_j$
is the matrix of orthonormal contrast coefficients for $\alpha_j$,
and $\varepsilon$ is the random error.\vadjust{\goodbreak}
Denote $n_j = (s-1)^j {{n}\choose{j}}$ and $X_j =
(x_{ik}^{(j)})_{N\times n_j}$,
then the (generalized) word-length pattern of design $D$ can be defined by
%
\begin{equation} \label{eq:A_j}
A_j(D) = N^{-2} \sum_{k=1}^{n_j}\Biggl|\sum_{i=1}^N
x_{ik}^{(j)}\Biggr|^2 \qquad\mbox{for } j=0,\ldots,n.
\end{equation}

For two designs $D^{(1)}$ and $D^{(2)}$, $D^{(1)}$ is
said to have less aberration than~$D^{(2)}$ if there exists an
$r\in\{1,2,\ldots,n\}$,
such that $A_r(D^{(1)})< A_r(D^{(2)})$ and\break $A_i(D^{(1)})= A_i(D^{(2)})$ for
$i=1,\ldots,r-1$. $D^{(1)}$ is said to have (generalized)
minimum aberration if there is no other design with
less aberration than~$D^{(1)}$.

For a regular design, the traditional definition of $A_j(D)$ is the
number of words of length $j$. Following \citet{XW2001}, $A_j(D)$
defined in (\ref{eq:A_j}) is the number of degrees of freedom
associated with all words of length $j$. Therefore, two definitions are
equivalent and generalized minimum aberration reduces to minimum
aberration for regular designs. For simplicity, in the following we use
the notion of word-length pattern and minimum aberration for both
regular and nonregular designs.

For an $(N, s^n)$-design $D=(x_{ik})_{N\times n}$, its centered
$L_2$-discrepancy (CD) is defined as
%
\begin{eqnarray}\label{cdf}
\qquad\phi(D) &=&\frac{1}{N^2} \sum^{N}_{i=1}\sum
^{N}_{j=1}\prod^n_{k=1}
\biggl(1+\frac{1}{2}\bigg|u_{ik}-\frac{1}{2}\bigg|+\frac{1}{2}
\bigg|u_{jk}-\frac{1}{2}\bigg|\nonumber
-\frac{1}{2}|u_{ik}-u_{jk}|\biggr)\\[-8pt]\\[-8pt]
&&{}-\frac{2}{N}{\sum^{N}_{i=1}\prod^n_{k=1}
\biggl(1+\frac{1}{2}\bigg|u_{ik}-\frac{1}{2}\bigg|-\frac{1}{2}
\bigg|u_{ik}-\frac{1}{2}\bigg|^2
\biggr)+\biggl(\frac{13}{12}\biggr)^n},\nonumber
\end{eqnarray}
where $u_{ik}=(2x_{ik}+1)/(2s)$. Note that $0< u_{ik} < 1$.

It is well known that word-length pattern
remains the same for combinatorially isomorphic designs.
However, the centered $L_2$-discrepancy will not be the same
when levels of factors are permuted.

\begin{example} \label{ex1}
Consider the two designs given in Table \ref{Tab_TwoDes}.
Both designs have one word of length three and share the
same word-length pattern $(A_1, A_2, A_3)=(0,0,2)$.
Their CD values are $0.033186$ and $0.033034$,
respectively. So Design $B$ is better than Design $A$ in terms of CD.
\end{example}

There is a close relationship between minimum aberration and uniformity
for two-level designs. \citet{FM2000} showed that
for a two-level regular design $D$, the
centered $L_2$-discrepancy of $D$ can be linearly expressed by
its word-length pattern $(A_1(D), A_2(D),\ldots,A_n(D))$.
Later on, \citet{MF2001} generalized it
to the nonregular case. Obviously, their results
cannot be generalized to high-level designs.

For an $s$-level factor, there are $s!$ possible level permutations.
Given an $(N,s^n)$-design $D$, we apply all $s!$ level permutations to
each column
and obtain $(s!)^n$ combinatorially isomorphic designs.
Denote the set of these designs as $\mathcal{P}(D)$.\vadjust{\goodbreak}
Some of them
{may be}
geometrically nonisomorphic and have different CD values.
We compute the CD value for each design and define $\bar{\phi}(D) $
as the average CD value of all designs in $\mathcal{P}(D)$, that is,
\[
\bar{\phi}(D) = \frac{1}{(s!)^n}\sum_{D' \in\mathcal{P}(D)}{\phi
(D')}.
\]

Note that all designs in $\mathcal{P}(D)$ share the same word-length
pattern. The following result shows that the average CD value, $\bar
{\phi}(D)$, is closely related to the word-length pattern of
$D$.\vspace*{-2pt}

\begin{theorem}\label{The_Relation}
For an $(N,3^n)$-design $D$,
\[
{\bar{\phi}(D)}
=\biggl(\frac{13}{12}\biggr)^n
- \biggl(\frac{29}{27}\biggr)^n+ \biggl(\frac{29}{27}\biggr)^n
\sum_{i=1}^n \biggl(\frac{2}{29}\biggr)^{i}A_i(D).\vspace*{-2pt}
\]
\end{theorem}

Theorem \ref{The_Relation} implies that the average
centered $L_2$-discrepancy and the minimum aberration criterion
are approximately equivalent,
as $(2/29)^i$ decreases geometrically when $i$ increases. Thus designs
permuted from a minimum aberration design tend to
be more likely to have low discrepancies. As will be seen, Theorem~\ref{The_Relation} is very useful in finding uniform FF designs.\vspace*{-2pt}

\begin{example}\label{ex:18}
Consider designs from the commonly used orthogonal array $\operatorname{OA}(18,3^7)$;
see, for example, Table 2(a) of \citet{XCW2004}. There are 3,
4, 4 and 3 combinatorially nonisomorphic designs when projected onto
3, 4, 5 and 6 factors, respectively. We rank these designs based upon
the minimum aberration criterion, and denote them as 18-3.1, 18-3.2,
18-3.3 and etc. For each design, we conduct all possible level
permutations and compute their CD values. Table~\ref{Tab_ACD18} shows
the average, minimum, maximum and standard deviation of the CD values
of all permuted designs, as well as one representative of the columns,
and $A_3$ and $A_4$ of the word-length pattern. Note that $A_1=A_2=0$
for all designs here. It can be seen from Table~\ref{Tab_ACD18} that
the rankings of average, minimum and maximum CD values are all
consistent with the minimum aberration ranking; that is, less
aberration leads to lower CD values. It is interesting to note that
designs 18-4.3 and \mbox{18-4.4} have the same word-length pattern but
different standard deviations; so do designs 18-6.2 and 18-6.3. This
implies that the word-length pattern does not uniquely determine the
variance of the CD values of permuted designs.\vspace*{-2pt}
\end{example}

We further compare the minimum aberration designs with the uniform
designs listed on the Uniform Design (UD) homepage
(\texttt{%
\href{http://www.math.hkbu.edu.hk/UniformDesign/}{http://www.math.hkbu.}
\href{http://www.math.hkbu.edu.hk/UniformDesign/}{edu.hk/UniformDesign/}}).
These uniform designs, labeled as UD18-3, UD18-4, etc., appear to be
orthogonal arrays of strength 2 so that $A_1=A_2=0$. The minimum
aberration design 18-3.1 has the same minimum CD value as the uniform
design UD18-3; however, the former has less aberration ($A_3=0.5$ vs.\vadjust{\goodbreak}
$A_3=0.67$) than the latter. Design 18-4.1  and UD18-4 have the same
properties, and they are indeed combinatorially isomorphic. Design
18-5.1 has a slightly larger minimum CD value (0.065265 vs. 0.065248)
and less aberration ($A_3=5$ vs. $A_3=6.17$) than UD18-5. The same
phenomenon also appears for design 18-6.1 and UD18-6. UD18-7 has a
smaller CD value than design 18-7 although they have the same
word-length pattern. The existing uniform designs have minimum
discrepancy for all cases because the run size is small here;
nevertheless, the level-permuted minimum aberration designs are
competitive. In summary, by permuting minimum aberration designs from
the commonly used $\operatorname{OA}(18,3^7)$, we can obtain good uniform FF designs.

\begin{table}
\caption{Comparison of $18$-run designs}\label{Tab_ACD18}
\begin{tabular*}{\textwidth}{@{\extracolsep{\fill}}lcccccd{2.2}d{2.2}@{}}
\hline
\textbf{Design} & \textbf{Columns} & \textbf{Ave} ${\bolds{\phi}}$ & \textbf{Min} ${\bolds{\phi}}$ & \textbf{Max} ${\bolds{\phi}}$ & \textbf{Sd} ${\bolds{\phi}}$
 & \multicolumn{1}{c}{$\bolds{A_3}$} & \multicolumn{1}{c@{}}{$\bolds{A_4}$} \\
\hline
18-3.1 & 1 2 3& 0.032526 & 0.032500 & 0.032538 & 0.000018 & 0.5 & \\
18-3.2 & 1 2 5& 0.032729 & 0.032500 & 0.032958 & 0.000163 &1 & \\
18-3.3 & 1 3 4& 0.033135 & 0.033034 & 0.033186 & 0.000072 &2& \\
UD18-3 & & & 0.032500 & & & 0.67& \\[3pt]
18-4.1 & 2 3 4 5& 0.047407 & 0.047357 & 0.047446 & 0.000023 &2& 1.5 \\
18-4.2 & 1 2 3 5& 0.047611 & 0.047391 & 0.047866 & 0.000166 & 2.5 & 1 \\
18-4.3 & 1 2 3 4& 0.048017 & 0.047849 & 0.048077 & 0.000087 & 3.5 & 0 \\
18-4.4 & 1 2 5 6& 0.048017 & 0.047849 & 0.048306 & 0.000139 &3.5 & 0 \\
UD18-4 & & & 0.047357 & & &2& 1.5\\[3pt]
18-5.1 & 2--6& 0.065273 & 0.065265 & 0.065337 & 0.000019 & 5 & 7.5 \\
18-5.2 & 1--3 5 6& 0.065883 & 0.065706 & 0.066193 & 0.000150 & 6.5 & 4.5 \\
18-5.3 & 1--5& 0.066086 & 0.065722 & 0.066423 & 0.000197 & 7 & 3.5 \\
18-5.4 & 1 2 5--7& 0.066492 & 0.066197 & 0.067107 & 0.000211 & 8 & 1.5\\
UD18-5 & & & 0.065248 & & &6.17& 5.17 \\[3pt]
18-6.1 & 2--7& 0.086964 & 0.086914 & 0.087145 & 0.000057 & 10 & 22.5 \\
18-6.2 & 1--6& 0.088184 & 0.087769 & 0.088591 & 0.000215 & 13 & 13.5 \\
18-6.3 & 1--3 5--7& 0.088184 & 0.087769 & 0.088974 & 0.000240 & 13 & 13.5 \\
UD18-6 & & & 0.086896 & & &12.33& 15.5 \\[3pt]
18-7 & 1--7 & 0.115386 & 0.114505 & 0.116556 & 0.000347 & 22 & 34.5 \\
UD18-7 & & & 0.113591 & & & 22 & 34.5 \\
\hline
\end{tabular*}
\vspace*{-3pt}
\end{table}

As suggested by Example \ref{ex:18}, an efficient way for constructing
uniform FF designs is to start with a minimum aberration design,
permute its
levels and choose the level permutation with the minimum CD value.
These designs have {minimum aberration, and good uniformity},
and are suitable for investigation of both nominal and quantitative factors.

\begin{definition}\label{Def_UMA}
Let $D$ be a minimum aberration design. If $D_* \in\mathcal{P}(D)$
has the minimum centered $L_2$-discrepancy over $\mathcal{P}(D)$, then
$D_*$ is said to be a uniform minimum aberration design.\vadjust{\goodbreak} 
\end{definition}

\begin{example}\label{ex:27}
Consider $27$-run designs. For $n=4$ to $13$ columns,
we evaluate average CD values for the existing
regular minimum aberration designs [see \citet{Xu2005}] and compare
with the CD values of the best designs listed on the UD homepage.
For $n=8$ to $10$, the average CD values of the minimum aberration designs
are even smaller than the CD values of the best existing designs; see
Table~\ref{Tab_ACD27} above.\vspace*{-3pt} 
\end{example}

To find uniform minimum aberration designs, we further conduct all possible
level permutations for these minimum aberration designs and calculate
the minimum and maximum CD values. Table \ref{Tab_ACD27} shows the
comparison between permuted minimum aberration (PMA) designs and the
best designs listed on UD homepage in terms of discrepancy and
aberration. For all designs, $A_1=0$ is not listed in the table.
For $n=4$, the PMA design is geometrically
isomorphic to the one listed on UD homepage. For $n=5$, the PMA design
has a larger CD value than the one listed on UD homepage, but the PMA
design has less aberration. For $n=6$, the PMA design has the same CD
value as the
one listed on UD homepage and has less aberration.
For $n> 6$, PMA designs always outperform the best
ones listed on UD homepage. Note that those designs listed on UD
homepage have resolution~2 ($A_2 >0$) whereas our designs have
resolution 3 ($A_2=0$), when $n > 6$.
This shows the advantage of our approach and the disadvantage
of the purely algorithmic approach. Further notice for $n=8$ and $9$,
even the maximum CD values of all permuted designs, are less
than those of the best existing ones.\looseness=-1\vspace*{-3pt}

\begin{table}
\caption{Comparison of $27$-run designs}\label{Tab_ACD27}
\begin{tabular*}{\textwidth}{@{\extracolsep{\fill}}d{2.0}d{1.7}d{1.7}d{1.7}ce{3.0}cd{2.2}d{2.2}@{}}
\hline
\multicolumn{6}{@{}c}{\textbf{Minimum aberration designs}} & \multicolumn{3}{c@{}}{\textbf{Designs on UD homepage}}\\[-5pt]
\multicolumn{6}{@{}c}{\hrulefill} & \multicolumn{3}{c@{}}{\hrulefill}\\
\multicolumn{1}{@{}l}{$\bolds{n}$} & \multicolumn{1}{c}{\textbf{Ave} ${\bolds{\phi}}$} & \multicolumn{1}{c}{\textbf{Min} ${\bolds{\phi}}$} & \multicolumn{1}{c}{\textbf{Max} ${\bolds{\phi}}$}
 & \multicolumn{1}{c}{$\bolds{A_2}$} & \multicolumn{1}{c}{$\bolds{A_3}$} &\multicolumn{1}{c}{${\bolds{\phi}}$} & \multicolumn{1}{c}{$\bolds{A_2}$} & \multicolumn{1}{c@{}}{$\bolds{A_3}$} \\
\hline
4 & 0.046549 & 0.046547\tabnoteref[\diamond]{tz} & 0.046553 & 0 & 0 & 0.046547 &0 & 0 \\
5 & 0.063818 & 0.063689& 0.063878 & 0 & 2 & 0.063525 & 0 & 2.67 \\
6 & 0.083786 & 0.083475\tabnoteref[\diamond]{tz} & 0.083923 & 0 & 4 & 0.083475 &0 & 5.33 \\
7 & 0.108701 & 0.108061\tabnoteref{dz} & 0.109118 & 0 & 10 & 0.108698 &0.10 & 12.17 \\
8 & 0.137749\tabnoteref{dz} & 0.136644\tabnoteref{dz} & 0.138483\tabnoteref{dz}& 0 & 16 & 0.138657 & 0.35 & 18.44 \\
9 & 0.172783\tabnoteref{dz} & 0.170996\tabnoteref{dz} & 0.174090\tabnoteref{dz}& 0 & 24 & 0.175343 & 0.69 & 30.05 \\
10 & 0.218927\tabnoteref{dz} & 0.213994\tabnoteref{dz} & 0.221241 & 0 & 42 &0.219131 & 1.36 & 40.99 \\
11 & 0.273255 & 0.264549\tabnoteref{dz} & 0.276195 & 0 & 60 & 0.272383 & 2& 56 \\
12 & 0.338698 & 0.325027\tabnoteref{dz} & 0.343084 & 0 & 80 & 0.336401 &2.32 & 75.46 \\
13 & 0.418900 & 0.397890\tabnoteref{dz} & 0.425576 & 0 & 104 & 0.414783 &3.53 & 96.20\\
\hline
\end{tabular*}
\tabnotetext[\diamond]{tz}{The same CD value as the best existing design;}
\tabnotetext[*]{dz}{Smaller CD value than the best existing design.}
\vspace*{-3pt}
\end{table}

\section{Construction of regular three-level uniform minimum
aberration designs}\label{Sec_ThreeLevel}
This section is devoted to providing an efficient method
for constructing uniform minimum aberration designs.\vadjust{\goodbreak}
For an $(N,3^n)$-design $D$, the total number of designs
in $\mathcal{P}(D)$ is $6^n$. However, when $D$ is a
regular FF design, many designs in $\mathcal{P}(D)$ are geometrically
isomorphic and have the same CD values. So it will be
much easier to find the uniform FF design
when a regular minimum aberration design is permuted.

For a three-level factor, exchange of
levels $0$ and $2$ does not change the geometrical
structure and such a ``mirror image'' operation keeps its centered
$L_2$-discrepancy unchanged according to formula (\ref{cdf}).
Denote $\pi_{i_0 i_1 i_2}$ as
a~permutation of $(0,1,2)$, that is, $\pi_{i_0 i_1 i_2}$ maps $(0, 1,
2)$ to $(i_0, i_1, i_2)$.
In view of the ``mirror image'' operation,
{we only
need to consider three permutations $\pi_{012}$, $\pi_{120}$ and $\pi
_{201}$ for a three-level design.}
Notice that $\pi_{012}$ is the identity map, $\pi_{120}$ maps $x$ to
$x+1 \pmod{3}$ and
$\pi_{201}$ maps $x$ to $x+2 \pmod{3}$.
So each permutation is equivalent to a linear permutation, which
transforms~$x$ to $x+b \pmod{3}$, where $b=0$, 1, or 2.

A regular $3^{n-k}$ FF design $D$ has $n-k$ independent columns,
denoted as $\mathbf{x}_1, \ldots, \mathbf{x}_{n-k}$, and $k$ dependent
columns, denoted as $\mathbf{y}_1, \ldots, \mathbf{y}_k$. It is specified by
$k$ linear equations:
\[
\cases{
\mathbf{y}_1  = c_{11}\mathbf{x}_1  +  c_{12}\mathbf{x}_2  +  \cdots + c_{1,n-k}\mathbf{x}_{n-k}  +b_1,\cr
\mathbf{y}_2  =  c_{21}\mathbf{x}_1  +  c_{22}\mathbf{x}_2  +  \cdots + c_{2,n-k}\mathbf{x}_{n-k}  +b_2,\cr
 \cdots  \cr
\mathbf{y}_k  =  c_{k1}\mathbf{x}_1  + c_{k2}\mathbf{x}_2  +  \cdots + c_{k,n-k}\mathbf{x}_{n-k}  +  b_k,}
\]
where $c_{ij}$ and $b_i$ are constants in $\operatorname{GF}(3)$, the finite field of
size 3.
Note that here and after, all algebra operations are performed
{in $\operatorname{GF}(3)$.}
{The standard design corresponds to $b_1=\cdots=b_k=0$ and} is an
$(n-k)$-dim linear space over $\operatorname{GF}(3)$.
Now any linear permutation of factor levels only alters
the coefficient vector $(b_1,\ldots,b_k)^T$.
Obviously, designs corresponding to the same
vector $(b_1,\ldots,b_k)^T$ are actually the same.
Thus among all the $3^n$ linearly permuted designs,
there are at most $3^k$ intrinsic differences. Moreover,
each design corresponding to
a specific $(b_1,\ldots,b_k)^T$
can be obtained by only conducting linear permutations
to the $k$ dependent columns (mapping $\mathbf{y}_j$ to $\mathbf{y}_j+b_j$
for $j=1,\ldots,k$), while keeping the
$n-k$ independent columns unchanged.
So we have the following lemma.

\begin{lemma}
For a regular $3^{n-k}$ FF design, {when all possible linear level
permutations are considered,
the set of all $3^n$ permuted designs consists of $3^{n-k}$ copies of
the $3^k$ designs obtained by permuting the $k$ dependent columns}.
\end{lemma}

For a design corresponding to vector $(b_1,\ldots,b_k)^T$, {consider
the ``mirror image'' permutation for all factors, that is, substituting
${\x}_i$ by $(2-\mathbf{x}_i)$ for $i=1,\ldots, n-k$ and $\mathbf{y}_j$ by
$(2-\mathbf{y}_j)$ for $j=1,\ldots, k$. The resulting ``mirror image''
design actually corresponds to the coefficient
vector $(2-2\sum_{i=1}^{n-k}c_{1i}-b_1, \ldots,2-2\sum
_{i=1}^{n-k}c_{ki}-b_{k})^T$. Because the ``mirror image'' permutation
does not change\vadjust{\goodbreak} the geometrical structures, these two designs are
geometrically isomorphic and have the same CD value. When two
coefficient vectors are the same, these two designs are identical.}
Thus we have the following lemma.

\begin{lemma}\label{Lem_GeoIso}
For a regular $3^{n-k}$ FF design, there are at most
{$(3^k+1)/2$}
geometrically nonisomorphic designs when all possible level
permutations are considered.
\end{lemma}

Applying the above results, we conduct level permutations of
three-level minimum aberration designs with 27 runs and 81 runs given
by \citet{Xu2005}
to find designs with minimum discrepancy. The results are concisely
presented as follows. For $27$-run designs,
when $n = 4$ to $6$, the first $n$ columns of
$\mathbf{x}_1$, $\mathbf{x}_2$, $\mathbf{x}_3$, $\mathbf{x}_1+\mathbf{x}_2+\mathbf{x}_3+2$,
$\mathbf{x}_1+2\mathbf{x}_2+1$ and $\mathbf{x}_1+\mathbf{x}_2+2\mathbf{x}_3+1$
form a regular uniform minimum aberration design;
when $n=7$ to $13$, the first $n$ columns of
$\mathbf{x}_1$, $\mathbf{x}_2$, $\mathbf{x}_3$, $\mathbf{x}_1+\mathbf{x}_2+\mathbf{x}_3+1$,
$\mathbf{x}_1+2\mathbf{x}_2+1$, $\mathbf{x}_1+\mathbf{x}_2+2\mathbf{x}_3$,
$\mathbf{x}_1+\mathbf{x}_3+2$, $\mathbf{x}_2+2\mathbf{x}_3+1$,
$\mathbf{x}_1+2\mathbf{
x}_2+2\mathbf{x}_3+2$,
$\mathbf{x}_1+\mathbf{x}_2+2$, $\mathbf{x}_2+\mathbf{x}_3+2$,
$\mathbf{x}_1+2\mathbf{
x}_2+\mathbf{x}_3$ and $\mathbf{x}_1+2\mathbf{x}_3+1$
form a regular uniform minimum aberration design,
where $\mathbf{x}_1, \mathbf{x}_2$, $\mathbf{x}_3$ are independent columns.
Their CD values are listed as ``Min $\phi$'' in
Table \ref{Tab_ACD27}.

For $81$-run designs, according to \citet{Xu2005},
when $n=5$ to $11$, the first
$n$ columns of
$\mathbf{x}_1$, $\mathbf{x}_2$, $\mathbf{x}_3$, $\mathbf{x}_4$, $\mathbf{ x}_1+\mathbf{
x}_2+\mathbf{ x}_3+\mathbf{ x}_4$,
$\mathbf{ x}_1+2\mathbf{ x}_2+\mathbf{ x}_3$, $\mathbf{ x}_1+2\mathbf{ x}_3+\mathbf{ x}_4$,
$\mathbf{ x}_1+2\mathbf{ x}_2+2\mathbf{ x}_4$, $\mathbf{ x}_2+\mathbf{ x}_3+2\mathbf{ x}_4$,
$\mathbf{ x}_1+\mathbf{ x}_2+2\mathbf{ x}_3+2\mathbf{ x}_4$ and $\mathbf{ x}_1+\mathbf{ x}_2$
form the minimum aberration design; when $n=12$ to $20$,
the first $n$ columns of
$\mathbf{ x}_1$, $\mathbf{ x}_2$, $\mathbf{ x}_3$, $\mathbf{ x}_4$, $\mathbf{ x}_1+\mathbf{
x}_2+\mathbf{ x}_3+\mathbf{ x}_4$,
$\mathbf{ x}_1+2\mathbf{ x}_2+\mathbf{ x}_3$, $\mathbf{ x}_1+2\mathbf{ x}_3+\mathbf{ x}_4$,
$\mathbf{ x}_1+2\mathbf{ x}_2+2\mathbf{ x}_4$, $\mathbf{ x}_1+\mathbf{ x}_2$, $\mathbf{
x}_2+2\mathbf{ x}_3+\mathbf{ x}_4$,
$\mathbf{ x}_1+2\mathbf{ x}_2+2\mathbf{ x}_3$, $\mathbf{ x}_1+2\mathbf{ x}_3+2\mathbf{ x}_4$,
$\mathbf{ x}_1+\mathbf{ x}_3$, $\mathbf{ x}_1+2\mathbf{ x}_2+\mathbf{ x}_4$, $\mathbf{
x_2}+\mathbf{ x}_3$,
$\mathbf{ x}_1+\mathbf{ x}_2+\mathbf{ x}_3+2\mathbf{ x}_4$, $\mathbf{ x}_1+\mathbf{ x}_2+2\mathbf{ x}_3$,
$\mathbf{ x}_2+2\mathbf{ x}_3+2\mathbf{ x}_4$, $\mathbf{ x}_1+\mathbf{ x}_4$ and $\mathbf{
x}_2+\mathbf{ x}_4$
form the minimum aberration design, where $\mathbf{ x}_1, \mathbf{ x}_2, \mathbf{
x}_3$ and $\mathbf{ x}_4$ are independent columns.
Table \ref{Tab_ACD81} summarizes the results when the minimum
aberration designs are permuted.
For example, when $n=7$, the best linear permutation conducted
to three dependent columns is $(0, 2,1)$, which means that
the best design with minimum CD value $0.102515$
is formed by seven columns
$\mathbf{ x}_1$, $\mathbf{ x}_2$, $\mathbf{ x}_3$, $\mathbf{ x}_4$, $\mathbf{ x}_1+\mathbf{
x}_2+\mathbf{ x}_3+\mathbf{ x}_4$,
$\mathbf{ x}_1+2\mathbf{ x}_2+\mathbf{ x}_3+2$ and $\mathbf{ x}_1+2\mathbf{ x}_3+\mathbf{ x}_4+1$.
As three-level designs with 81 runs are not listed on UD homepage,
the best designs found in Table \ref{Tab_ACD81} are apparently new.

\begin{table}
\tabcolsep=5.5pt
\caption{{Results of $81$-run minimum aberration designs}}\label{Tab_ACD81}
\begin{tabular*}{\textwidth}{@{\extracolsep{\fill}}lccl@{}}
\hline
$\bolds{n}$ &\textbf{Ave} ${\bolds{\phi}}$ & \textbf{Min} $\bolds{\phi}$ &\multicolumn{1}{c@{}}{\textbf{Best level permutations}} \\
\hline
\phantom{0}5& {0.062691} & {0.062690}&
\begin{tabular}{cccccccccc@{}}
0 & & & & & & & & & \\
\end{tabular}
\\
\phantom{0}6
& {0.081294} & {0.081290}&
\begin{tabular}{cccccccccc@{}}
0 & 1 & & & & & & & & \\
\end{tabular}
\\
\phantom{0}7
& {0.102528} & {0.102515} &
\begin{tabular}{cccccccccc@{}}
0 & 2 & 1 & & & & & & & \\
\end{tabular}
\\
\phantom{0}8
& {0.126795} & {0.126764} &
\begin{tabular}{cccccccccc@{}}
0 & 2 & 1 & 0 & & & & & & \\
\end{tabular}
\\
\phantom{0}9
& {0.154565} & {0.154497} &
\begin{tabular}{cccccccccc@{}}
0 & 2 & 1 & 0 & 1 & & & & & \\
\end{tabular}
\\
10
& {0.186393} & {0.186255} &
\begin{tabular}{cccccccccc@{}}
0 & 2 & 1 & 0 & 1 & 0 & & & & \\
\end{tabular}
\\
11
& {0.226648} & {0.225969} &
\begin{tabular}{cccccccccc@{}}
1 & 1 & 0 & 0 & 0 & 0 & 2 & & & \\
\end{tabular}
\\
12
& {0.270884} & {0.269750} &
\begin{tabular}{cccccccccccccccc@{}}
1 & 1 & 0 & 0 & 2 & 0 & 0 & 2 & & & &&&&&\\
\end{tabular}
\\
13
& {0.324370} & {0.322305} &
\begin{tabular}{cccccccccccccccc@{}}
1 & 0 & 2 & 0 & 2 & 0 & 2 & 1 & 2 & &&&&&& \\
\end{tabular}
\\
14
& {0.385994} & {0.382976} &
\begin{tabular}{cccccccccccccccc@{}}
0 & 0 & 1 & 1 & 2 & 2 & 2 & 0 & 2 & 2 &&&&&&\\
\end{tabular}
\\
15
& {0.457704} & {0.453338} &
\begin{tabular}{cccccccccccccccc@{}}
0 & 0 & 1 & 1 & 2 & 2 & 2 & 0 & 2 & 2 & 2 &&&&& \\
\end{tabular}
\\
16
& {0.540883} & {0.534813} &
\begin{tabular}{cccccccccccccccc@{}}
0 & 0 & 1 & 1 & 2 & 2 & 2 & 0 & 2 & 2 & 2 & 2 & &&& \\
\end{tabular}
\\
17
& {0.640085} & {0.631437} &
\begin{tabular}{cccccccccccccccc@{}}
0 & 0 & 1 & 1 & 2 & 2 & 2 & 0 & 2 & 2 & 2 & 2 & 0 &&& \\
\end{tabular}
\\
18
& {0.755854} & {0.743782} &
\begin{tabular}{cccccccccccccccc@{}}
0 & 0 & 1 & 1 & 2 & 2 & 2 & 0 & 2 & 2 & 2 & 2 & 0 & 1 & & \\
\end{tabular}
\\
19
& {0.898270} & {0.883749} &
\begin{tabular}{cccccccccccccccc@{}}
0 & 0 & 1 & 1 & 2 & 2 & 2 & 0 & 2 & 2 & 2 & 2 & 0 & 1 & 2 & \\
\end{tabular}
\\
20
& {1.066298} & {1.048120} &
\begin{tabular}{cccccccccccccccc@{}}
0 & 0 & 1 & 1 & 2 & 2 & 2 & 0 & 2 & 2 & 2 & 2 & 0 & 1 & 2 & 1 \\
\end{tabular}
\\
\hline
\end{tabular*}
\vspace*{-3pt}
\end{table}

As stated in Lemma \ref{Lem_GeoIso},
for a regular $3^{n-k}$ FF design, there are at most
{$(3^k+1)/2$}
geometrically nonisomorphic designs when all possible level
permutations are considered.
Now consider the simplest case with $k=1$.
A regular $3^{n-1}$ minimum aberration design can be specified
by $\y=2\x_1+\cdots+2\x_{n-1}+b$, where $\x_1, \ldots, \x_{n-1}$
are the independent columns, $\y$ is the dependent column and $b \in
\operatorname{GF}(3)$. This is equivalent to $\x_1+\cdots+\x_{n-1}+\y=b \pmod{3}$.
When levels of $\y$ are permuted, they will generate
two geometrically nonisomorphic designs.
To be specific, denote $D_{i}$ as the design corresponds to $b = n+i
\pmod{3}$ for $i=0,1,2$.
Then $D_0$ contains a row of ones, and $D_{1}$ and $D_{2}$ are
geometrically isomorphic.
The next theorem provides explicit formulas for the CD values of
$D_{0}$ and $D_{1}$.

\begin{theorem}\label{The_3^n-11}
Let $D_{0}$ and $D_{1}$ be the two geometrically nonisomorphic
regular $3^{n-1}$ minimum aberration designs, where $D_{0}$ represents
the design with the all-one row.
Then the centered $L_2$-discrepancies of $D_{0}$ and $D_{1}$ are
\[
\phi(D_{0}) = \biggl(\frac{13}{12}\biggr)^n
- \biggl(\frac{29}{27}\biggr)^n+2 \biggl(\frac{2}{27}\biggr)^n
+\frac{2(-1)^{n}}{3^{3n}}
\]
and
\[
\phi(D_{1}) = \biggl(\frac{13}{12}\biggr)^n
- \biggl(\frac{29}{27}\biggr)^n+ 2\biggl(\frac{2}{27}\biggr)^n +
\frac{(-1)^{n+1}}{3^{3n}}.
\]
\end{theorem}

As an immediate implication of Theorem \ref{The_3^n-11},
when $n$ is odd, $\phi(D_{1})>\phi(D_{0})$, so $D_{0}$
is the uniform minimum aberration design;
when $n$ is even, $\phi(D_{0})>\phi(D_{1})$, so $D_{1}$
is the uniform minimum aberration design.

\section{\texorpdfstring{Connection to the $\beta$-word-length pattern}{Connection to the beta-word-length pattern}}\label{Sec_Justify}
Under the hierarchical principle, \citet{CY2004} defined the
$\beta$-word-length pattern. Specifically,
for an $(N,s^n)$-design, let
$\theta_j$ be the vector of all $j$-degree interactions and~$Z_j$
be the matrix of orthogonal polynomial contrast coefficients for
$\theta_j$.
Then the response $Y$ can be fitted by a polynomial model
${Y} = Z_0 \theta_0 + Z_1 \theta_1 + \cdots+ Z_K \theta_K +
\varepsilon$.
Denote $Z_j = (z_{ik}^{(j)})_{N\times n'_j}$, where $n'_j$ is the
number of effects with degree $j$. The $\beta$-word-length pattern is
defined by
\[
\beta_j(D) = N^{-2} \sum_{k=1}^{n'_j}\Biggl|\sum_{i=1}^N
z_{ik}^{(j)}\Biggr|^2 \qquad\mbox{for } j=0,\ldots,K,\vadjust{\goodbreak}
\]
where $K=n(s-1)$ represents the highest polynomial degree. \citet{CY2004} argued that a good design should minimize $\beta_1, \beta_2,
\ldots, \beta_K$ in a~sequential order.

It is interesting to see how uniform minimum aberration designs perform
under the $\beta$-word-length pattern.
In principle, given a design,
one can always find the best design related to the $\beta$-word-length
pattern by
permuting levels for all factors.
However, the computational burden
makes it infeasible to evaluate all $\beta_j(D)$ values
even for three-level designs with moderate number of factors.
Here we only consider permuting levels of regular minimum aberration
designs with $27$ runs and $n=4$ to $10$ columns and compute their
$\beta$-word-length patterns.
To our surprise, for all cases, the permuted designs with best
$\beta$-word-length patterns always have the least centered $L_2$-discrepancies
and vice verse; that is, the uniform minimum aberration designs are the
best designs
under the $\beta$-word-length pattern.
Of course, there are cases where different designs have the same
CD value but different $\beta$-word-length pattern, and vice verse.
Moreover, for $n=4$ to $8$, the $\beta$-word-length pattern
and centered $L_2$-discrepancy give the exactly same ordering of the
permuted designs.
For $n=9$ or 10, the orderings under the two criteria are quite consistent,
though not identical.

We end this section with a theoretical result. Notice that a regular
$3^{n-1}$ minimum aberration design has resolution $n$ so that
$A_1=\cdots=A_{n-1}=0$, which implies $\beta_1=\cdots=\beta_{n-1}=0$.
The following theorem gives an interesting relationship between the CD
value and $\beta_n$.

\begin{theorem}\label{The_3^n-1.2}
For a regular $3^{n-1}$ minimum aberration design $D$,
\[
\phi(D) = \biggl(\frac{13}{12}\biggr)^n
- \biggl(\frac{29}{27}\biggr)^n+ 2 \biggl(\frac{2}{27}\biggr)^n - 2
\biggl(\frac{1}{27}\biggr)^n + \biggl(\frac{2}{27}\biggr)^n \beta_n(D).
\]
\end{theorem}

Theorem \ref{The_3^n-1.2} shows that the two criteria, centered
$L_2$-discrepancy and $\beta$-word-length pattern, are exactly
equivalent for regular $3^{n-1}$ minimum aberration designs.

\section{Conclusion}\label{con}

Uniform FF designs are useful for studying quantitative factors with
multiple levels; however, the construction of such designs is
challenging. We establish a connection between uniformity and
aberration by showing that the average centered $L_2$-discrepancy is a
function of the word-length pattern. We propose to construct uniform FF
designs by permuting levels of existing minimum aberration designs.
Using this strategy, we construct regular uniform minimum aberration
designs with 27 runs and 81 runs for practical use. We further evaluate
the performance of the uniform minimum aberration designs for
polynomial models. They perform well under the $\beta$-word-length pattern.

\begin{appendix}\label{app}
\section*{Appendix: Proofs of all theorems}

For an $(N,s^n)$-design $D=(x_{ik})$, let
$d_H(i,j)$ be the Hamming distance of rows $i$ and $j$ of $D$, that is,
$d_H(i,j)=\sharp\{k \dvtx x_{ik}\not=x_{jk}, k=1, \ldots, n\}$, where
$\sharp(S)$ is the cardinality of $S$.
The distance distribution of $D$ is $(B_0(D),B_1(D),\ldots,\allowbreak B_n(D))$, where
\[
B_j(D) = N^{-1} \sharp\{(a,b)\dvtx d_H(a,b) =j \mbox{ and } a, b=1,\ldots
,N\}\qquad
\mbox{for } j=0, \ldots, n.
\]

\citet{XW2001} showed that the (generalized) word-length pattern can
be calculated by
the MacWilliams transform of the distance distribution, that is,
\[
A_j(D)=N^{-1} \sum_{i=0}^n B_{i}(D) P_j(i;n,s)\qquad
\mbox{for } j=0,\ldots,n,
\]
where $P_j(x;n,s) = \sum_{i=0}^j (-1)^i(s-1)^{j-i}
{{x}\choose{i}}{{n-x}\choose{j-i}}$ are the Krawtchouk polynomials.

By the orthogonality of the Krawtchouk polynomials, we also have
\[
B_j(D)=N\cdot s^{-n}\sum_{i=0}^n P_j(i;n,s) A_i(D).
\]

The following existing property related to
Krawtchouk polynomials was stated in \citet{MS1977}.
\begin{lemma}\label{Lem_Krawtchouk}
For nonnegative integers $n, x$ and $s$ with $n \ge x$, $s \ge2$ and
$0<y<1$,
\[
\sum_{j=0}^n P_j(x;n,s) y^j = [1+(s-1)y]^{n-x}(1-y)^x.
\]
\end{lemma}

To prove Theorem \ref{The_Relation}, we need the following lemma.

\begin{lemma}\label{Lem_ReCD2-WD}
For an $(N,s^n)$-design $D$, denote $\delta_{ij}$ as the number of
positions where rows $i$ and $j$
take the same value, that is, $\delta_{ij}=n-d_H(i,j)$.
Then for any positive number $z$ greater than 1,
\[
\sum_{i,j=1}^N z^{\delta_{ij}} =N^2 \biggl(\frac{z+s-1}{s}
\biggr)^{n} \sum_{i=0}^n \biggl(\frac{z-1}{z+s-1}\biggr)^i A_i(D).
\]
\end{lemma}

\begin{pf}
According to the definition of distance distribution,
Lemma \ref{Lem_Krawtchouk} and the relationship between distance
distribution and word-length pattern, we have
\begin{eqnarray*}
 \sum_{i,j=1}^N z^{\delta_{ij}} & = &
N \sum_{j=0}^n B_j(D)z^{n-j}=z^n s^{-n}N^2\sum_{j=0}^n \sum_{i=0}^n
P_j(i;n,s) A_i(D)z^{-j} \\
& = &  {z^n s^{-n}N^2 \sum_{i=0}^n
\biggl(1+\frac{s-1}{z}\biggr)^{n-i}\biggl(1-\frac{1}{z}\biggr)^i
A_i(D)}\\
& = &  N^2 \biggl(\frac{z+s-1}{s}\biggr)^{n} \sum
_{i=0}^n \biggl(\frac{z-1}{z+s-1}\biggr)^i A_i(D).
\end{eqnarray*}
\upqed
\end{pf}

\begin{pf*}{Proof of Theorem \ref{The_Relation}}
Notice that for an $(N,3^n)$-design $D=(x_{ik})$ with $x_{ik}=0,1$, or
2, $u_{ik}$ and $u_{jk}$ in
formula (\ref{cdf}) can only take values $1/6$, $1/2$ and~$5/6$.
If $u_{ik} = 1/2$, $1+\frac{1}{2}|u_{ik}-\frac{1}{2}|
-\frac{1}{2}|u_{ik}-\frac{1}{2}| ^{2}$ takes value $1$;
otherwise, it takes value $10/9$. Furthermore,
if $u_{ik} = u_{jk} = 1/6$ or $u_{ik} = u_{jk} = 5/6$,
$1+\frac{1}{2}|
u_{ik}-\frac{1}{2}|+\frac{1}{2}| u_{jk}-\frac{1}{2}| -\frac{1}{2}|
u_{ik}-u_{jk}|$ takes value $4/3$; otherwise, it takes value~$1$.
Thus for any two rows of an $(N,3^n)$-design $D$, denoted by
$(x_{i1},x_{i2},\ldots,x_{in})$ and
$(x_{j1},x_{j2},\ldots,x_{jn})$, if we define
$\gamma_i(D)= \sharp\{k \dvtx x_{ik} \not= 1, k=1,\ldots,n
\}$ and $\gamma_{ij}(D)= \sharp\{k \dvtx x_{ik}=x_{jk} \not=1
,k=1,\ldots,n\}$,
the CD value of $D$ can be determined by
the distributions of its $\gamma_i$ and $\gamma_{ij}$ values. That is,
formula (\ref{cdf}) can be simplified to
%
\begin{eqnarray}\label{Eq_1}
\phi(D)&=&{\biggl(\frac{13}{12}\biggr)}^n-\frac{2}{N}
{\sum_{i=1}^{N}
{\biggl(\frac{10}{9}\biggr)}^{\gamma_i(D)}}\nonumber\\[-8pt]\\[-8pt]
&&{}+\frac
{1}{N^2}{ {\sum_{i=1}^{N}{\biggl(\frac{4}{3}
\biggr)}^{\gamma_i(D)}}}
+\frac{1}{N^2}{ {\sum_{i\not=j}{\biggl(\frac
{4}{3}\biggr)}^{\gamma_{ij}(D)}}}.\nonumber
\end{eqnarray}

Moreover, for any fixed row $i$ of $D$, when all level permutations
of $D$ are considered, each $n$-tuple in $\operatorname{GF}(3)^n$ occurs $2^n$ times.
For any element $x_{ik}$ in the $k$th column, there are three possible
choices, that is,
$0,1,2$, corresponding to $\frac{10}{9}, 1, \frac{10}{9}$ for
$1+\frac{1}{2}|u_{ik}-\frac{1}{2}|
-\frac{1}{2}|u_{ik}-\frac{1}{2}| ^{2}$.
So
\[
{\sum_{D' \in\mathcal{P}(D)}{\biggl(\frac{10}{9}\biggr)}^{\gamma
_i(D')} =
2^n \cdot\biggl(\frac{10}{9} +\frac{10}{9} + 1\biggr)^n = 2^n\cdot
{\biggl(\frac{29}{9}\biggr)}^n}.
\]
Similarly,
\[
{\sum_{D' \in\mathcal{P}(D)}{\biggl(\frac{4}{3}\biggr)}^{\gamma
_i(D')} =
2^n \cdot\biggl(\frac{4}{3} +\frac{4}{3} + 1\biggr)^n = 2^n\cdot
{\biggl(\frac{11}{3}\biggr)}^n=6^n\cdot{\biggl(\frac{11}{9}\biggr)}^n}.
\]
For any two rows $i$ and $j$ of $D$ with
Hamming distance $d_H(i,j)=n-\delta_{ij}$, when level permutations
of corresponding columns are considered, each identical pair, that is,
$(l,l), l\in \operatorname{GF}(3)$,
occurs twice in $\delta_{ij}$ positions where rows $i$ and $j$
coincide, and each distinct pair
occurs once in corresponding $d_H(i,j)$ positions where rows $i$ and
$j$ differ.
So
\[
{\sum_{D' \in\mathcal{P}(D)}{\biggl(\frac{4}{3}\biggr)}^{\gamma_{ij}(D')}}=
2^{\delta_{ij}}\cdot\biggl(\frac{4}{3} +\frac{4}{3} + 1
\biggr)^{\delta_{ij}} \cdot(6\times1)^{n-\delta_{ij}}
= 6^n \biggl(\frac{11}{9}\biggr)^{\delta_{ij}}.
\]
Then
\begin{eqnarray*}
 {\sum_{D'\in\mathcal{P}(D)}\phi(D')}&=&
 {6^n \biggl(\frac{13}{12}\biggr)^n
- \frac{2}{N} \cdot N\cdot2^n \biggl(\frac{29}{9}\biggr)^n + \frac
{1}{N^2} 6^n \sum_{i,j=1}^N \biggl(\frac{11}{9}\biggr)^{\delta
_{ij}}} \\
&=& {6^n \biggl(\frac{13}{12}\biggr)^n - 2\cdot6^n
\biggl(\frac{29}{27}\biggr)^n + \frac{6^n}{N^2} \sum_{i,j=1}^N
\biggl(\frac{11}{9}\biggr)^{\delta_{ij}}}
\end{eqnarray*}
and
\[
 {\bar{\phi}(D)}=
 {\biggl(\frac{13}{12}\biggr)^n - 2 \biggl(\frac
{29}{27}\biggr)^n +
\frac{1}{N^2} \sum_{i,j=1}^N\biggl(\frac{11}{9}\biggr)^{\delta_{ij}}}.
\]
Finally, the result follows from Lemma \ref{Lem_ReCD2-WD} and the fact $A_0=1$.
\end{pf*}

\begin{pf*}{Proof of Theorem \ref{The_3^n-11}}
We will use formula (\ref{Eq_1}) to calculate the CD values.
First, we prove that the distributions of
$\gamma_{ij}$ values for the two designs $D_{0}$ and $D_{1}$
are the same. Because $D_{1}$ is obtained by adding $1 \pmod{3}$ to
the last column of $D_{0}$, $\gamma_{ij}(D_{0}) \ne\gamma
_{ij}(D_{1})$ if and only if both last positions
of rows $i$ and $j$ of $D_{0}$ have the same value $0$ or $1$.
For any two distinct rows of $D_{0}$ with the last positions
both taking value $0$, denoted by
$(x_{i1},x_{i2},\ldots,x_{i,n-1},0)$ and $(x_{j1},x_{j2},\ldots
,x_{j,n-1},0)$, respectively,
there exists a unique pair of rows of $D_{1}$,
$(x_{i1},\ldots, x_{i,l-1},x_{il}-1,x_{i,l+1},\ldots,x_{i,n-1},2)$
and $(x_{j1},\ldots, x_{j,l-1},x_{jl}-1, x_{j,l+1},\ldots,x_{j,n-1},2)$,
where $l$ is the first position such that $x_{il} \not= x_{jl}$.
These two pairs have the same $\gamma_{ij}$ value.
Similarly, for any two distinct rows of $D_{0}$ with the last positions
both taking value $1$,
$(x_{i1},x_{i2},\ldots,x_{i,n-1},1)$
and $(x_{j1},x_{j2},\ldots,x_{j,n-1},1)$,
there exists a unique pair of
rows of $D_{1}$, $(x_{i1},\ldots, x_{i,l-1}, x_{il}+1,x_{i,l+1},\ldots
,x_{i,n-1},1)$
and $(x_{j1},\ldots, x_{j,l-1},$ $x_{jl}+1, x_{j,l+1},\ldots
,x_{j,n-1},1)$, with
the same $\gamma_{ij}$ value, where $l$ is the
first position such that $x_{il} \not= x_{jl}$.
It is easy to see that the above correspondence between two pairs of
rows in $D_{0}$ and $D_{1}$ is actually one-to-one (bijective). This
completes our claim on the distributions of
$\gamma_{ij}$ values.

Now we consider the distributions of $\gamma_i$ values of $D_{0}$ and $D_{1}$.
For a design~$D$, denote $\eta_j(D) = \sharp\{i \dvtx \gamma_i(D) = j,
i=1, \ldots, N\}, j=0,1,\ldots,n$,
as the distribution of $\gamma_i$ values of $D$. Obviously, $\eta
_0(D_{0}) = 1$ as
$D_{0}$ contains the all-one row, and
$\eta_0(D_{1}) = 0$. Moreover, for any
$j=0,1,\ldots,n$, the total number of rows in $\operatorname{GF}(3)^n$ with
$\gamma_i = j$ is ${n \choose j} 2^{j}$; therefore,
$\eta_j(D_{0})+2\eta_j(D_{1}) = {n \choose j} 2^{j}$.
For convenience, if a vector with length $j$
only contains $0$ or $2$, it will be called a
$(0,2)^j$-vector.
Then ${\eta_j(D_{0})}/{{n \choose j}}$ is the number of
possible $(0,2)^j$-vectors with
sum $j \pmod{3}$ and ${\eta_j(D_{1})}/{{n \choose j}}$ is the
number of
possible $(0,2)^j$-vectors with
sum {$j+1 \pmod{3}$}. Notice that a $(0,2)^j$-vector with
sum {$j-1 \pmod{3}$} can be obtained
by conducting a ``mirror image'' operation to
a $(0,2)^j$-vector with
sum {$j+1 \pmod{3}$}.
Thus
${\eta_j(D_{1})}/{{n \choose j}}$ also represents the number of
possible $(0,2)^j$-vectors with
sum {$j-1 \pmod{3}$}.
Each $(0,2)^j$-vector with sum $j \pmod{3}$
can be formed by adding $2$ or $0$ to a $(0,2)^{j-1}$-vector
with sum $(j-1)-1$ or $(j-1)+1 \pmod{3}$. So we have
${\eta_j(D_{0})}/{{n \choose j}} = 2 {\eta_{j-1}(D_{1})}/{{n \choose j-1}}$.
Combining this with $\eta_j(D_{0})+2\eta_j(D_{1}) = {n \choose j} 2^{j}$
and $\eta_0(D_{1}) = 0$, we obtain $\eta_j(D_{1}) = {n \choose
j}\frac{2^j-(-1)^j}{3}$ and
$\eta_j(D_{0}) = {n \choose j}\frac{2^j+2(-1)^j}{3}$. Thus $\eta
_j(D_{0}) - \eta_j(D_{1}) ={n \choose j} (-1)^j$.
Using formula (\ref{Eq_1}) and $N=3^{n-1}$, we have
\begin{eqnarray*}
\phi(D_{0})-\phi(D_{1})
&=&  {\frac{1}{N^2}\sum_{j=0}^n \pmatrix{n \cr j}(-1)^j
\biggl(\frac{4}{3}\biggr)^j
- \frac{2}{N}\sum_{j=0}^n \pmatrix{n \cr j}(-1)^j \biggl(\frac
{10}{9}\biggr)^j } \\
&=&  {\frac{1}{N^2}\biggl[\biggl(-\frac{1}{3}\biggr)^n
- 2N \biggl(-\frac{1}{9}\biggr)^n \biggr] } \\
&=& {\frac{3^n - 2N}{N^2}\biggl(-\frac{1}{9}\biggr)^n =
\frac{(-1)^n}{3^{3n-1}}}.
\end{eqnarray*}
From Theorem \ref{The_Relation}, we also have
\[
\phi(D_{0}) + 2 \phi(D_{1}) = 3\biggl[\biggl(\frac{13}{12}\biggr)^n
- \biggl(\frac{29}{27}\biggr)^n+ 2\biggl(\frac{2}{27}
\biggr)^n\biggr].
\]
Then we have
\[
\phi(D_{0}) = \biggl(\frac{13}{12}\biggr)^n
- \biggl(\frac{29}{27}\biggr)^n+2 \biggl(\frac{2}{27}\biggr)^n
+\frac{2(-1)^{n}}{3^{3n}}
\]
and
\[
\phi(D_{1}) = \biggl(\frac{13}{12}\biggr)^n
- \biggl(\frac{29}{27}\biggr)^n+ 2\biggl(\frac{2}{27}\biggr)^n +
\frac{(-1)^{n+1}}{3^{3n}}.
\]
\upqed
\end{pf*}

\begin{pf*}{Proof of Theorem \ref{The_3^n-1.2}}
 For a regular $3^{n-1}$
minimum aberration design $D=(x_{ik})$ with resolution $n$, its $\beta
_n(D)$ is determined by the product of linear polynomials as follows:
%
\begin{equation}\label{eq:betan}
\beta_n(D) = N^{-2} \Biggl|\sum_{i=1}^N p_1(x_{i1}) \times\cdots
\times p_1(x_{in})\Biggr|^2 ,
\end{equation}
where $N=3^{n-1}$ and $p_1(x)=\sqrt{3/2}(x-1)$.
Because $p_1(x)=0$ when $x=1$, we only need to consider rows with 0 or
2 only, that is, $(0,2)^n$-vectors.
Notice that $D_{0}$ is an
$(n-1)$-dim linear space or a coset over $\operatorname{GF}(3)$ containing the all-one vector,
thus run $(2-z_1,\ldots, 2-z_n)$ occurs in $D_{0}$ if and only
if $(z_1,\ldots,z_n)$ occurs in~$D_{0}$.
So for any odd $n$, $\beta_n(D_{0})=0$ according to (\ref{eq:betan}).
To calculate $\beta_n(D_{1})$ for odd~$n$, we will establish
a recursive formula.
For this purpose, we use
$D_i^n$ to denote a design $D_i$ with $n$ columns; that is,
the sum of each row of the design is congruent to $n+i$ modulo 3 for $i=0,1,2$.
Then, up to row permutations, we can express $D_1^n$ as follows:
%
\begin{equation}\label{eq:D_1^n}
D_{1}^n = \left[\matrix{
D_{2}^{n-1} & \zero\cr
D_1^{n-1} & \one\cr
D_0^{n-1} & \two\cr
}
\right].
\end{equation}
Let $\delta{(D_i^{n-1})}$ be the difference between the number of
$(0,2)^{n-1}$-vectors in $D_i^{n-1}$ with even number of zeros and the
number of those with odd number of zeros for $i=0,1,2$. Then, according
to (\ref{eq:betan}), we have
%
\begin{equation}\label{eq:betan1}
\beta_{n-1}(D_{i}^{n-1})=3^{-2(n-2)} (3/2)^{n-1} | \delta
{(D_i^{n-1})} |^2 \qquad\mbox{for } i=0,1,2.
\end{equation}
Furthermore, when $n-1$ is even, $\delta{(D_0^{n-1})} + \delta
{(D_1^{n-1})} + \delta{(D_2^{n-1})} =0$ and $\delta{(D_0^{n-1})}
=-2\delta{(D_1^{n-1})} = -2\delta{(D_2^{n-1})}$.
Then, according to (\ref{eq:betan}) and (\ref{eq:D_1^n}), for odd
$n$, we have
%
\begin{eqnarray}\label{eq:betan2}
\beta_n(D_{1}^n) &=& 3^{-2(n-1)} (3/2)^n |- \delta(D_{2}^{n-1}) +
\delta(D_0^{n-1}) |^2 \nonumber
\\[-8pt]\\[-8pt]
&=& 3^{-2(n-1)} (3/2)^n | -3 \delta(D_{1}^{n-1}) |^2 .\nonumber
\end{eqnarray}
Combining (\ref{eq:betan1}) and (\ref{eq:betan2}), we obtain $ \beta
_n(D_{1}^n) = (3/2) \beta_{n-1}(D_{1}^{n-1}).$
In the same vein, for even $n$, we can establish
$\beta_n(D_{1}^n) = (1/6)\beta_{n-1}(D_1^{n-1})$
and $\beta_n(D_{0}^n) = (2/3)\beta_{n-1}(D_1^{n-1})$.
{So for odd $n$, $\beta_n(D_{1}^n) = (3/2) \beta_{n-1}(D_{1}^{n-1}) =
(1/4)\times  \beta_{n-2}(D_1^{n-2})$.}
It is easy to verify that $\beta_3(D_1^{3}) = 3/8$; thus for odd $n$,
$\beta_n(D_{1}^n) = (3/8)(1/4)^{(n-3)/2}= 3/2^n$.
For even $n$, we obtain $\beta_n(D_{1}^n) = (1/6)(3/2^{n-1}) =
1/{2^n}$ and $\beta_n(D_{0}^n) =(2/3)(3/2^{n-1})=1/2^{n-2}$.
Then the result follows from Theorem~\ref{The_3^n-11}.
\end{pf*}
\end{appendix}

\section*{Acknowledgments}
This research was done when the first
author was visiting the Department of Statistics at University of
California, Los Angeles.
{The authors would like to thank an associate editor and a referee for
their comments and suggestions which helped to improve this paper.}

%

\printaddresses

\end{document}